\documentclass[9 pt]{article}
\usepackage{amssymb}
\usepackage{amsfonts}
\usepackage{amsmath}
\usepackage{graphicx}

\providecommand{\U}[1]{\protect\rule{.1in}{.1in}}
\numberwithin{equation}{section}
\numberwithin{equation}{section}
\numberwithin{equation}{section}
\def\e{\varepsilon}

\newtheorem{theorem}{Theorem}[section]

\newtheorem{definition}[theorem]{Definition}

\newtheorem{lemma}[theorem]{Lemma}

\newtheorem{proposition}[theorem]{Proposition}

\newtheorem{remark}[theorem]{Remark}

\newenvironment{proof}[1][Proof]{\noindent\textbf{#1.} }{\ \rule{0.5em}{0.5em}}

\begin{document}

\title{Integral representation results in $BV\times L^p$.} 
\author{\textsc{Gra\c{c}a Carita}\thanks{%
CIMA-UE, Departamento de Matem\'{a}tica, Universidade de \'{E}vora, Rua Rom%
\~{a}o Ramalho, 59 7000 671 \'{E}vora, Portugal e-mail: gcarita@uevora.pt}, 
\textsc{\ Elvira Zappale}\thanks{%
D.I.In., Universita' degli Studi di Salerno, Via Giovanni Paolo II 132,
84084 Fisciano (SA) Italy e-mail:ezappale@unisa.it}}
\maketitle

\begin{abstract}
An integral representation result is obtained for the relaxation of a class
of energy functionals depending on two vector fields with different
behaviors, which may appear in the context of image decomposition and termochemical equilibrium problems. 


\end{abstract}

\section{Introduction}

Minimization of energies depending on two independent vector fields have been introduced to model several phenomena. Namely, when $u$ is a Sobolev function  in $W^{1,q}, q >1,$ and $v$ is in $L^p$, the study of these energies (see \eqref{functional}) was motivated by the analysis of coherent thermochemical equilibria in a multiphase multicomponent system, with $\nabla u$ representing the elastic strain and $v$ the chemical composition of the material. In the theory of linear magnetostriction, the stored energy depends on the linearized strain and the direction of magnetization, we refer to \cite{FKP2, FKP1} and the references therein for more details. Moreover, when $p=q$ this type of energies is used to model Cosserat theory and bending phenomena in nonlinear  elasticity and also for the description of thin structures, see \cite{LDR} and \cite{BFM}. Here $v$ takes into account either Cosserat vectors or bending moments and $\nabla u$ is the elastic strain. When $u$ is a function of bounded variation, functionals similar to \eqref{functional} enter into image decomposition models, i.e., in order to denoise and restore a given image $f$, it is required to minimize a functional which is the sum of a 'total variation' term (i.e. a norm of $Du$) and a penalization term, i.e. a norm in a suitable functional space of $f-u-v$. Essentially $f$ can be decomposed into the sum of two components $u$ and $v$.  The first component (cartoon), $u$, is well structured and it describes the homogeneous objects which are present in the image. The second component, $v$, contains the oscillating pattern (both textures and noise), and we refer to \cite{ROF, Meyer, VO,VO2, AABFC} among the wide literature in this field.

In order to cover a wide class of applications we start from the functional setting $W^{1,1}\times L^p$, with anysotropic energies with linear growth in the gradient variable $\nabla u$. Indeed, let $1<p\leq\infty$, for every $(u,v) \in W^{1,1}( \Omega;%
\mathbb{R}^{n}) \times L^{p}(\Omega;\mathbb{R}^{m}) $ 
define the functional 
\begin{equation}  \label{functional}
J( u,v) :=\int_{\Omega}f(x,u(x), v(x),\nabla u(x)) dx
\end{equation}
where $f:\Omega\times\mathbb{R}^{n}\times\mathbb{R}^{m}\times\mathbb{R}%
^{n\times N}\rightarrow[ 0,+\infty) $ is a continuous function. We discuss separately the cases $1<p<\infty$ and $p=\infty.$ Thus we introduce for $1<p<\infty$ the functional%
\begin{equation}
\overline{J}_{p}( u,v) :=\inf\{ \underset{n\rightarrow \infty%
}{\lim\inf}J(u_{n},v_{n}) :u_{n}\in W^{1,1}( \Omega;\mathbb{%
R}^{n}) ,~v_{n}\in L^{p}(\Omega;\mathbb{R}^{m})
,~u_{n}\rightarrow u\text{ in }L^{1},~v_{n}\rightharpoonup v\text{ in }%
L^{p}\} ,   \label{relaxedp}
\end{equation}
for any pair $(u,v) \in BV( \Omega;\mathbb{R}^{n})
\times L^{p}( \Omega;\mathbb{R}^{m}) $ and, for $p=\infty$ the
functional%
\begin{equation}
\overline{J}_{\infty}( u,v) :=\inf\{ \underset{n\rightarrow
\infty}{\lim\inf}J(u_{n},v_{n}) :u_{n}\in W^{1,1}( \Omega;%
\mathbb{R}^{n}) ,~v_{n}\in L^{\infty}( \Omega ;\mathbb{R}%
^{m}) ,~u_{n}\rightarrow u\text{ in }L^{1},~v_{n}\overset{\ast}{%
\rightharpoonup}v\text{ in }L^{\infty}\} ,   \label{relaxedinfty}
\end{equation}
for any pair $(u,v) \in BV(\Omega;\mathbb{R}^{n})
\times L^{\infty}(\Omega;\mathbb{R}^{m}).$

Since bounded sequences $\{u_n\}$ in $W^{1,1}(\Omega;\mathbb{R}^m)$ converge
in $L^1$ to a $BV$ function $u$ and bounded sequences $\{v_n\}$ in $%
L^p(\Omega;\mathbb{R}^m)$ if $1<p\leq \infty$, 
weakly converge to a function $v \in L^p(\Omega;\mathbb{%
R}^m)$, (weakly $\ast$ in $L^\infty$), the relaxed functionals $\overline{J}%
_{p}$ and $\overline{J}_{\infty}$ will be composed by a Lebesgue part, a
jump part concentrated on the jump set of $u\in BV(\Omega;\mathbb{R}^m)$ and
a Cantor part, absolutely continuous with respect to the Cantor part of the
distributional gradient $Du$. On the other hand, as already emphasized in 
\cite{FKP2}, it is crucial to observe that $v$ is not defined on the jump
and the 'Cantor part' sets of $u$, thus specific features of the density $f$
will come into play to ensure a proper integral representation. The one of (\ref{relaxedp}) is obtained in \cite{CZ}, via the blow-up method introduced in \cite{FM2},
under the following hypotheses:

\noindent $(H_0)$ $f(x,u,\cdot, \cdot)$ is convex-quasiconvex for every $%
(x,u)\in \Omega \times \mathbb{R}^n$;

\begin{itemize}
\item[$(H_1)_{p}$] There exists a positive constant $C$ such that%
\begin{equation*}
\frac{1}{C}( \vert b\vert ^{p}+\vert \xi\vert
) -C\leq f( x,u,b,\xi) \leq C( 1+\vert
b\vert ^{p}+\vert \xi\vert) \hbox{ for every }(x,u,b,\xi) \in\Omega\times\mathbb{R}^{n}\times \mathbb{R}%
^{m}\times\mathbb{R}^{n\times N};
\end{equation*}

\item[$( H_{2}) _{p}$] For every compact set $K\subset\Omega
\times\mathbb{R}^{n}$ there exists a continuous function $\omega _{K}:%
\mathbb{R\rightarrow}[ 0,+\infty) $ with $\omega_{K}(
0) =0$ such that

\begin{itemize}
 \item[$(1)$] $\vert f( x,u,b,\xi) -f(
x^{\prime},u^{\prime},b,\xi) \vert \leq\omega_{K}(\vert x-x^{\prime }\vert +\vert u-u^{\prime}\vert) ( 1+\vert b\vert ^{p}+\vert \xi\vert) $
for every $(x,u,b,\xi) $ and $( x^{\prime},u^{\prime
},b,\xi) \in K\times\mathbb{R}^{m}\times\mathbb{R}^{n\times N};$
\item[$(2)$] Moreover, given $x_{0}\in\Omega$ and $\varepsilon>0$ there exists $\delta>0$
such that if $\vert x-x_{0}\vert <\delta$ then 
\begin{equation*}
\displaystyle{f(
x,u,b,\xi) -f( x_{0},u,b,\xi) \geq-\varepsilon(
1+\vert b\vert ^{p}+\vert \xi\vert),} \hbox{ for every }(u,b,\xi) \in\mathbb{R}^{n}\times\mathbb{R}^{m}\times%
\mathbb{R}^{n\times N};
\end{equation*} 

\end{itemize}
\item[$( H_{3})_{p}$] There exist $c^{\prime}>0,~L>0,~0<\tau%
\leq1$, $0<s<p$ such that%
\begin{equation*}
t>0,~\xi\in\mathbb{R}^{n\times N},\text{ with }t\vert \xi\vert
+ t|b|^p>L\Longrightarrow\left\vert \frac{f( x,u,t^{\frac{1}{p}}b,t\xi) }{t}-f_p^{\infty
}( x,u,b,\xi) \right\vert \leq c^{\prime}\left( \frac{\vert
b\vert ^{s}}{t^{1-\frac{s}{p}}} +\frac{\vert \xi\vert ^{1-\tau}}{t^{\tau}}%
\right) 
\end{equation*}
for every $(x,u) \in\Omega\times\mathbb{R}^{n},$ where $f_p^{\infty }$ is the $(p,1)$- recession function of $f$ defined for every $%
(x,u,b,\xi)\in \Omega \times \mathbb{R}^{n}\times \mathbb{R}^{m}\times \mathbb R^{N\times n}$ as 
\begin{equation}
\displaystyle{f_p^{\infty }(x,u,b,\xi ):=\limsup_{t\rightarrow +\infty}\frac{%
f(x,u,t^{\frac{1}{p}}b,t\xi )}{t}}.  \label{recessionp}
\end{equation}
\end{itemize}

\begin{theorem}
\label{MainResultp} Let $J$ and ${\overline J}_p$ be given by \eqref{functional} and  \eqref{relaxedp} respectively, with $f$ satisfying $(H_0)$, $(H_1)_p-(H_3)_p$ then%
$$
\overline{J}_p(u,v)= \int_\Omega f(x,u,v, \nabla u) dx +
\int_{J_u} K_p(x,0, u^+, u^-,\nu_u)d\mathcal{H}^{N-1} +
\int_\Omega f_p^\infty\left(x,u,0,\frac{dDu}{d|D^c u|}\right) d|D^c u|,
$$
for every $(u,v)\in BV(\Omega;\mathbb{R}^n)\times L^p(\Omega;\mathbb{R}^m)$,
where $K_p:\Omega \times\mathbb{R}%
^{n}\times\mathbb{R}^{n}\times\mathbb{R}^{m}\times S^{N-1}\rightarrow[
0,+\infty) $ is defined as
\begin{equation}
	K_p( x,c,d,b,\nu) :=\inf\left\{ \int_{Q_{\nu}}f_p^{\infty}(x,w( y),\eta( y) ,\nabla w(y))
	dy:w\in\mathcal{A}( c,d,\nu), \eta\in L^{\infty}( Q_{\nu};
	\mathbb{R}^{m}), \int_{Q_{\nu}}\eta(y)dy=b\right\}, 
	\label{Kp}
\end{equation}
with 
$$
\begin{array}{ll}
{\mathcal A}(c,d,\nu):=\{w \in W^{1,1}(Q_\nu;\mathbb R^n): &w(y)= c \hbox{ if } y\cdot \nu = \frac{1}{2}, w(y)= d  \hbox{ if } y \cdot \nu=-\frac{1}{2},\\
\\
&w \hbox{ is  1-periodic in }  \nu_1,\dots, \nu_{N-1} \hbox{ directions}\}.
\end{array}
$$
\end{theorem}

\noindent In order to provide an integral description of the functional $\overline{J}_{\infty }$, introduced
in $( \ref{relaxedinfty}) $, in \cite{CZ} we proved Theorem \ref{MainResultinfty} replacing assumptions $(
H_{1}) _{p}-( H_{3}) _{p}$ by the\ following ones:

\begin{itemize}
\item[$( H_{1}) _{\infty}$] Given $M>0,$ there exists $C_{M}>0$ 
such that, if $\vert b\vert \leq M$ then%
\begin{equation*}
\frac{1}{C_{M}}\vert \xi\vert -C_{M}\leq f(x,u,b,\xi)
\leq C_{M}( 1+\vert \xi\vert) , \hbox{ for every }(x,u,\xi) \in\Omega\times\mathbb{R}^{n}\times\mathbb{R%
}^{n\times N};
\end{equation*}

\item[$(H_2)_{\infty}$] For every $M>0,$ and for every compact set $%
K\subset\Omega\times\mathbb{R}^{n}$ there exists a continuous function 

$
\omega_{M,K}( 0) =0$ such that if $\vert b\vert \leq M$
then%
\begin{equation*}
\vert f(x,u,b,\xi) -f(
x^{\prime},u^{\prime},b,\xi) \vert \leq\omega_{M,K}(
\vert x-x^{\prime }\vert +\vert u-u^{\prime}\vert
) ( 1+\vert \xi\vert) 
\end{equation*}
for every $(x,u,\xi) , (x^{\prime},u^{\prime},\xi)
\in K\times\mathbb{R}^{n\times N}.$  Moreover, given $M>0,$ $x_{0}\in\Omega,$ and $\varepsilon>0$ there exists $%
\delta>0$ such that if $\vert b\vert \leq M$ and $\vert
x-x_{0}\vert \leq\delta$ then%
\begin{equation*}
f( x,u,b,\xi) -f( x_{0},u,b,\xi) \geq-\varepsilon
( 1+\vert \xi\vert)  \hbox{ for every }( u,\xi) \in\mathbb{R}^{n}\times\mathbb{R}^{n\times N};
\end{equation*}

\item[$(H_{3})_{\infty }$] Given $M>0$, there exist $c_{M}^{\prime
}>0,~L>0,~0<\tau \leq 1$ such that 
\begin{equation*}
\vert b\vert \leq M,~t>0,~\xi \in \mathbb{R}^{n\times N},~\text{%
with }t\vert \xi \vert >L\Longrightarrow \left\vert \frac{f(
x,u,b,t\xi) }{t}-f^{\infty }( x,u,b,\xi) \right\vert \leq
c_{M}^{^{\prime }}\frac{\vert \xi \vert ^{1-\tau }}{t^{\tau }}
\end{equation*}%
for every $( x,u) \in \Omega \times \mathbb{R}^{n},$
where $f^\infty(b,\xi)$  is the $(p,\infty)-$recession function, i.e. the `standard' recession function in the last variable, defined for every $%
(x,u,b,\xi)\in \Omega \times \mathbb{R}^{n}\times \mathbb{R}^{m}\times \mathbb R^{N\times n}$ as 
\begin{equation}
\label{recession}
\displaystyle{f^\infty(x,u,b,\xi):=\limsup_{t \to +\infty}\frac{f(x,u,b,t\xi)}{t}.}
\end{equation}
\end{itemize}

\begin{theorem}
\label{MainResultinfty} Let $J$ and ${\overline J}_\infty$ be given by \eqref{functional} and \eqref{relaxedinfty} respectively, with $f$
satisfying $(H_0)$, $(H_1)_{\infty}-(H_3)_{\infty}$ then
$$
\overline{J}_\infty(u,v)= \int_\Omega f(x,u,v,\nabla u) dx +
\int_{J_u} K_\infty(x,0, u^+, u^-,\nu_u)d\mathcal{H}^{N-1} +
\int_\Omega f^\infty\left(x,u,0,\frac{dDu}{d|D^c u|}\right) d|D^c u|,
$$
for every $(u,v)\in BV(\Omega;\mathbb{R}^n)\times L^p(\Omega;\mathbb{R}^m)$,
where $K_\infty:\Omega \times\mathbb{R}%
^{n}\times\mathbb{R}^{n}\times\mathbb{R}^{m}\times S^{N-1}\rightarrow[
0,+\infty) $  is defined by 
\begin{equation}
K_\infty( x,c,d,b,\nu) :=\inf\left\{ \int_{Q_{\nu}}f^{\infty}(x,w( y),\eta( y) ,\nabla w(y))
dy:w\in\mathcal{A}( c,d,\nu), \eta\in L^{\infty}( Q_{\nu};
\mathbb{R}^{m}), \int_{Q_{\nu}}\eta(y)dy=b\right\}. 
\label{Kinfty}
\end{equation}
\end{theorem}

\section{Conclusions}

The above results generalize those contained in \cite{FKP2, FKP1} to the case where $f$ depends also on $x$ and $u$ and extend those contained in \cite{ADM, FM2} to the case where it appears in the energy an extra field $v$, which can be interpreted both as the chemical composition of a material, as a Cosserat director in elasticity theories, and as the noise (and texture) component in the imaging models (see for instance \cite{Meyer, ROF, VO, VO2}).  Moreover, Theorem \ref{MainResultp}  recovers \cite[Theorem 1.1]{RZCh}, where additive models, arising from imaging applications, are considered.

It is worth to observe that the functional $J$ in \eqref{functional} is well defined also for $v \in L^1(\Omega;\mathbb R^m)$ and the asymptotic behaviour of $\overline{J}_p$ can be studied when $p=1$ and the the weak topology in $L^p$ is replaced by the weak $\ast$ convergence in the sense of measures, thus leading to a relaxed functional defined in $BV(\Omega;\mathbb R^n)\times {\cal M}(\Omega;\mathbb R^m)$. Indeed, an integral representation can be deduced by \cite{BZZ} when the density $f$ does not depend on $x$ and $u$ and has linear growth in both variables $(b,\xi)$. In this latter case the integral representation is composed by four terms, being necessary to consider an integral with respect the part of the measure $v$ which might be singular with respect to $Du$ and the involved density contains the standard recession function $f^\infty$, defined as
$$
\displaystyle{f^\infty(b,\xi):=\limsup_{t \to +\infty}\frac{f(tb,t\xi)}{t}.}
$$

The proofs of Theorems \ref{MainResultp} and \ref{MainResultinfty} are contained in \cite{CZ} in the easier case when $f\equiv f(b,\xi)$. In this latter case both $K_\infty$ and $K_p$ reduce to the recession function  $f^\infty$ evaluated at $0$, namely 
$$
K_p(0, u^+, u^-, \nu_u)=K_\infty(0, u^+, u^-, \nu_u)= f^\infty\left(0,D^s u\right)= f^\infty(0,(u^+-u^-)\otimes \nu_u).
$$
Observe that, besides the set of assumptions on $f$, $(H_1)_p- (H_3)_p$, differs from $(H_1)_\infty -(H_3)_\infty$ and in general $f^\infty\not = f^\infty_p$, the equality between \eqref{recessionp} and \eqref{recession} holds when $b=0$.
We stress the fact that in \cite{CZ} other assumptions on $f$, which allow to achieve both theorems above, are also discussed. 
 
It is worth to observe also that assumption $(H_0)$ can be removed, thus replacing $f$ by its  convex-quasiconvex envelope in the above integral representations, \eqref{Kp} and \eqref{Kinfty}. 
\

\noindent{\bf Acknowledgements}

\noindent The authors thank Irene Fonseca and Ana Margarida Ribeiro for many helpful discussions on the topic of this article and the CNA at CMU for its kind hospitality. 
The research of the authors has been
partially supported by Funda{\c c}${\tilde {\rm a}}$o para a Ci${\hat {\rm e}}$ncia e Tecnologia (Portuguese Foundation) through CIMA-UE and GNAMPA-INdAM through project Un approccio variazionale all'analisi di modelli competitivi non lineari.

\end{document}